\documentclass[11pt]{article}

 \pdfoutput=1

\usepackage{latexsym}
\usepackage{amsfonts}
\usepackage{graphicx}
\usepackage{amssymb}
\usepackage{amsmath}
\usepackage{amscd}

\usepackage{amssymb,stmaryrd,tikz-cd}
\usetikzlibrary{
				matrix,
				arrows,	
				decorations.pathreplacing,	
				decorations.markings
}

	\pgfarrowsdeclarealias{<}{>}{to}{to}
	\pgfarrowsdeclaredouble{<<}{>>}{to}{to}
	\pgfarrowsdeclaretriple{<<<}{>>>}{to}{to}
	
	\tikzset{->-/.style={decoration={
	  markings,
	  mark=at position .5 with {\arrow{>}}},postaction={decorate}}}
	\tikzset{-->-/.style={decoration={
	  markings,
	  mark=at position .75 with {\arrow{>}}},postaction={decorate}}}
	\tikzset{->--/.style={decoration={
	  markings,
	  mark=at position .25 with {\arrow{>}}},postaction={decorate}}}
	\tikzset{->>-/.style={decoration={
	  markings,
	  mark=at position .55 with {\arrow{>>}}},postaction={decorate}}}
	\tikzset{->>>-/.style={decoration={
	  markings,
	  mark=at position .55 with {\arrow{>>>}}},postaction={decorate}}}

\usepackage[hidelinks=true]{hyperref}
\hypersetup{
  colorlinks  = true,   
  urlcolor    = black,  
  linkcolor   = black,   
  citecolor   = black,  
}
\date{  }

\title{Paper}

\author{Gabriel Zapata}

\newtheorem{Theorem}{Theorem}

\newtheorem{lemma}{Lemma}
\newtheorem{definition}{Definition}
\newtheorem{proposition}{Proposition}

\newcommand{\aom }
{
	\ensuremath
	{ 
	\phantom{ \mspace{2mu} }
			^{ \underline{ \phantom{ \mspace{10mu} } } }
	\phantom{ \mspace{2mu} }
	}
}

\newcommand{\sub }[1]
{
	\ensuremath
	{ 
			_{ \mspace{.5mu} _{#1}  }
	}
}

\newcommand{\zero}
{
\ensuremath
	{
		  _{ \mathrm{o} }  
	}
}

\newcommand{\ket}[1]
{
	\ensuremath
	{ 
		\langle
			\mspace{1mu} #1 \mspace{1mu}
		\rangle
	}
}

\newcommand{\by}
{
	\ensuremath
	{ 
		\!\cdot\!
	}
}

\newcommand{\sset}[1]
{ 
\ensuremath
	{ 
		 \{\mspace{.1mu} #1 \mspace{.2mu}\}
	}
}

\newcommand{\cset}[2]
{ 
\ensuremath
	{ 
		\{\, #1~|~ #2 \,\}
	}
}

\newcommand{\cat}[1]{
	\ensuremath{ 
		\mathsf{#1}
	}
}

\newcommand{\Set}{\mathsf{Set}}

\newcommand{\f}[1]{
	\ensuremath{ \mathrm{#1} }
}

\newcommand{\inv}[1]{
	\ensuremath{ {#1}^{^{\!\!\textrm{-}1}}\mspace{-5mu} }
}

\newcommand{\inverse}[0]{
	\ensuremath{ ^{^{-1} }\mspace{-4mu} }
}

\newcommand{ \kernel }[1]{
	\ensuremath{ 
		\mathrm{Ker}\, #1
	}
}

\newcommand{\im}[1]{
	\ensuremath{ 
		\mathrm{Im}\,#1 
	}
}

\newcommand{\mar}[2]{
	\ensuremath{ 
		#1\,
			\mspace{0.8mu}
				\mathbin{
					\tikz[baseline]
						\draw[->] 
						(0pt,.5ex) -- (3.2ex,.5ex);
				}
			\mspace{1.2mu}
		\,#2
	}
}

\newcommand{\marto}[2]{
	\ensuremath{ 
		#1\,
			\mspace{1.5mu}
				\mathbin{
					\tikz[baseline]
						\draw[|->] 
						(0pt,.5ex) -- (3.5ex,.5ex);
				}
			\mspace{2mu}
		\,#2
	}
}

\newcommand{\id}{
		\ensuremath{ 
 			 1\mspace{-4.7mu}\mathrm{l} 
	}
}

\newcommand{\idi}[1]{
\ensuremath{ 
		\mspace{1mu}
 			1\mspace{-4.6mu}\text{l}_{_{#1}} 
	}
}

\newcommand{\mor}[3]{
	\ensuremath{ \mathrm{Hom}  _{\, #1} (#2,#3) 
	}
}

\newcommand{\every}[2]
{
	\mbox
	{
					\ensuremath{\pmb{\forall}} 
						\ensuremath{\!_{#1\,\,\in\,#2\!}}
	}
}

\newcommand{\one}[2]
{
	\mbox
	{
				\ensuremath{\pmb{\exists}}
						\ensuremath{
						\mspace{-5mu}! _{\,\,#1\,\,\in\,#2\!}
						}
	}
}

\def\classification{
\@ifnextchar [{\@xfootnotenext}%
{\begingroup\let\protect\noexpand \xdef\@thefnmark{}\endgroup\@footnotetext}
}

\begin{document}

\title{Rewriting Group Products with Transversals}

\maketitle

\begin{abstract}
\noindent
For any  group $G$ with subgroup $H$ and a set of representatives $T$
from the set of cosets $G/H$, we develop  a  rewriting system from $G$ 
that bequeaths a product into the set decomposition  
 $T\times H$ of $G$, converting it into a group.
In return the rewriting system  also describes any product between 
elements in $G$ in terms of elements from $T$ and $H$, while we gain a new
universal group-isomorphism between  $G$ and $T\times H$. From this framework
we  develop the concept of Diffracted Groups.
\end{abstract}

\section{Introduction}

Let $H$ be a subgroup of a group $G$ and $T$ be a set of  
\emph{representatives} from the set of cosets
$G/H$. The set $T$, under the product of $G$,  is generally
a semigroup and not a group.
Therefore, our objective is to develop a rewriting process in $G$ 
bequeathing a group product for the set  $T\times H$, which is
 bijective to $G$, cf. 
 \cite{GZ0} chapter III. 
From the rewriting process we also demonstrate how the bijection between
$T\times H$ and $G$ lifts to a universal group-isomorphism.

\vspace{3pt}

To  formally design the rewriting process, we  begin by defining  
$\textsf{Set}$ and $\textsf{Grp}$ 
to be the  small category of  sets and groups, respectively.
Given any objects $H$ and $G$ from 
$\mathsf{C}\in \sset{\textsf{Set},\,\textsf{Grp}}$, 
we write
$\mor{\mathsf{C}}{H}{G}$ 
to denote the set of morphisms from $H$ to $G$.
We also let
$\f{S}\sub{H}$ be the group of permutation (or, the symmetric group)
 with
underlying set $H$. We can then express the  well-known, 
natural set-isomorphism
  \vspace{-4pt}
\begin{equation}\label{eq permutation adjunction}
\mor{\mathsf{Grp}}{G}{\f{S}\sub{H}} \,
	\;\cong\;
	  \mor{\mathsf{Set}}{G\times H}{H}\,
\end{equation}
\vspace{-18pt}

\noindent
between \emph{group-representations} 
 and  \emph{group-actions} on $H$.
This bridge will allow us to transverse between $G$'s local "representative" structure and its global "bundle" action on  $H$. Moreover, while also
 using the basic but 
 powerful \emph{Cayley's  representation},
we will go from the
internal multiplicative  structure of $G$ to the set-theoretic decomposition
$T\times H$ with respect to a subgroup $H$ and its
 representatives $T$ of $G/H$. 
We will then  develop a
representation of $G$  bequeathing a multiplicative rewriting process for
group elements in the  decomposition $T\times H$:

\begin{definition}\label{def  Cayley's representation}
\emph{The \emph{Cayley's  representation} is the monomorphism
defined as}
 \vspace{-3pt}
\begin{equation}\label{eq  cayleys representation}
\varrho\,:\,\mar{G}{\f{S}\sub{G} }
~~\mathrm{such ~that}~~
\every{ {g,\,h}\,}{\, G}
	~\big[~~	 
\varrho\,:\,
	\marto{g}{
		\big(~
			g^{\,\varrho}:\;\mar{h} {g\mspace{1mu}h}
		~\big) 
	}
~~\big]\,,
\end{equation}
\emph{where the element $g^{\,\varrho}\in\f{S}\sub{G}$ symbolizes the image of $g$ under $\varrho$.}
\end{definition}

\vspace{3pt}

To be more precise, Cayley's representation 
$\varrho\in \mor{\mathsf{Grp}}{G}{\f{S}\sub{G}} $  shows us a
decomposition of the product of $G$ by way of its associated (left) action 
$\hat{\varrho }\in \mor{\mathsf{Set}}{G\times G}{G}$, i.e., via the 
 commutative diagram
 \vspace{-3pt} 
\begin{equation}\label{dia  Cayley Product}
\begin{aligned}[c]
\begin{tikzpicture}[description/.style={fill=white}]
\matrix (m) [matrix of math nodes, row sep=1.8em, column sep=3em, 
			text height=1.5ex, text depth=0.25ex] 
	{ G \times G & G  \\
	    {G} \times G & G	\,,	\\ }; 
\path[->,font=\scriptsize]
	(m-1-1) edge node[left]  {$  \idi{_G} \times\idi{_G}$} (m-2-1)
	(m-1-1) edge node[above] {$ ``\bullet" $} 	(m-1-2)
    (m-2-1)edge node[auto]   {$ \hat{\varrho } $} (m-2-2);
     \path[<-,font=\scriptsize]
     	(m-1-2)edge node[right]  {$ \idi{_G} $} (m-2-2);
\end{tikzpicture}
\end{aligned}
\vspace{-5pt} 
\end{equation}
where ``$\bullet$" and  $\id\sub{G}$ denote  the  multiplicative 
operation and the identity automorphism of $ G $, respectively.
In the following sections, we develop representations that engender
the replacement of the pair ${G} \times G$, from Cayley's action 
$\hat{\varrho}$, by the pair
$
{(T\times H)} \times \,(T\times H) \,.
$ 
From the latter decomposition we will obtain a rewriting process describing
the product of elements from $G$ in terms of elements from $T$ and $H$. 
This rewriting process will also endow a multiplicative group structure 
for the set $T\times H$, and it will induce a
universal group-isomorphism  between $G$ and the group structure 
 developed for $T\times H$.

\section{The Frobenius Representation}
\label{sec  representatives diffraction}

This section reviews a method that offers an intuitive path for 
describing a  decomposition of $G$---a codification using a subgroup $H$ and 
$G/H$.

\vspace{-3pt}

\begin{definition}
\emph{Given a group $G$ with $H\leq G$ and its set of left cosets $G/H$, let}
\begin{equation}\label{def  power by quotient by}
\Lambda
	~:=~ 
		\{\,  
			\tau : \mar {G/H} 
				{
					\bigcup_{\, g  \,\in\; G}\,\,
					(gH)  
				}
			~|~
			\every{g } {\,G}~ [\; \tau(gH)\in gH \;]
		~\}\,.
\end{equation}
\emph{A set of $($left\,$)$ \emph{representatives}  $T$ of   $G/H$ 
 is the image of a coordinate $\tau\in \Lambda$, i.e.,
$
T \, : = \,\im \tau\,.
$
Moreover, a set of representatives $T$  is termed a 
\emph{transversal}  if  the identity element 
 $1\in G$ is also in $T$.}
\end{definition}

\vspace{-5pt}

Intuitively, a set of representatives  $T$ 
consists  of one chosen element 
from each coset of  $G/H$. 
The didactic definition is not an accident, it is suitable
for visualizing and
extrapolating crucial representative properties: 
 each coordinate point $\tau\in\Lambda$ contains a codification of $G$, 
modulo $H$. Given $\tau$ and the canonical projection 
$
p :\mar{G}{G/H}$, there is a unique map
$
\bar{\tau}:
 \mar{G}{T}
$ 
such that the following diagrams 
\begin{equation*}\label{dia  representatives} 
\begin{aligned}[c]
\begin{tikzpicture}[description/.style={fill=white}]
\matrix (m) [matrix of math nodes, row sep=.1 em, column sep=1.2em, 
			text height=1.5ex, text depth=0.25ex] 
	{  G    		&   & T  		   	   \\
	  \phantom{\!}  &   & \phantom{a}\mspace{60mu}     	    \\
	  		 		&   & G/H             	   \\
	}; 
\path[->, dotted, font=\scriptsize]
	(m-1-1)  edge node[above] {$\bar{\tau} $} (m-1-3);
\path[->,font=\scriptsize]
	(m-1-1)  edge node[below] {$ p $}      (m-3-3)	
	(m-3-3) edge node[right]  {$ \tau $}   (m-1-3);
\end{tikzpicture}
\end{aligned}
\text{,~i.e.,}\quad
\begin{aligned}[c]
\begin{tikzpicture}[description/.style={fill=white}]
\matrix (m) [matrix of math nodes, row sep=.1 em, column sep=1.2em, 
			text height=1.5ex, text depth=0.25ex] 
	{  g   		&   & 
						\tau(\,g\,H\,)\\
	  \phantom{\!}  &   & \phantom{\!}     	    \\
	  			&   &     g\,H            	  \\
	}; 
\path[|->, dotted, font=\scriptsize]
	(m-1-1)  edge node[above] {$  \bar{\tau}  $} (m-1-3);
\path[|->,font=\scriptsize]
	(m-1-1.300)  edge node[below] {$ p  $} (m-3-3.170)	
	(m-3-3) edge node[right] {$ \tau $} (m-1-3);
\end{tikzpicture}
\end{aligned}
\end{equation*}
commute, for any $g\in G$. In particular
$
 \bar{\tau}: \mar{G}{T}
$
systematically maps  $G$ onto $T$ (cf. \cite{Magnus} classical
 construction.)
\begin{definition}\label{def  representatives}
\emph{Let $p:\mar{G}{G/H}$ be a canonical projection and $\tau\in \Lambda$. 
Given its set of representatives  $T$ of $G/N$, a \emph{representative map}  
is the map}
\[
\overline{\tau}:\mar{G}{T}
\quad\mathrm{such ~that} \quad
\every{g}{G}
	\big[~
\begin{tikzcd}[column sep=26]
		  g   
				 \arrow[ |->,"\bar{\tau}" ]{r}  & 	
		  \overline{g}\,:=\, {\tau}(\,gH\,)
\end{tikzcd}
	~\big]\,.
\,
\]
\emph{Moreover any element of the form  $\overline{g}\in T$
is termed  a \emph{representative of $g$} (from $gH$.)}
\end{definition}

\begin{lemma}\label{lem  representative calculus}
If $T$ is a  set of representatives  of $G/H$ of $G$, then 
\begin{description}

\item[\cat{i\,.\phantom{i}}] 
$
\every{g}{G} 
	~\big[~~
		g\,H\,=\, \overline{g}\,H
	~~\big]
$

\item[\cat{ii\,.}]
$
\every{g\sub{1},\,g\sub{2}}{G} 
	~\big[~~
		 \overline{g\sub{1}\overline{g\sub{2}}} 
			\, = \,
				\overline{g\sub{1}\,g\sub{2}}
	~~\big]
$

\end{description}
\end{lemma}

\noindent\emph{Proof}: Let $T$ be a set of representative of
$G/H$ and    $\tau$ a map of its representatives.

\vspace*{4pt}

$\textsf{i.}$  Given $g\in G$, then $\, gH \, = \,\bar{g}\,H$ since
its representative $\bar{g}$  satisfies
\[
g\in \bar{g}\,H
	~\Longleftrightarrow~
		\one{h}{H}~[~\bar{g} \, = \,g\,h^{-1}~]
	~\Longleftrightarrow~		
		 \,\bar{g} \,\in gH \,.
\]

$\textsf{ii}.$ 
By the above argument
$
	\, g\sub{2} H\,=\, \overline{g\sub{2}}H
$ for any $g\sub{2}\in G$. Then
$
	 \,g\sub{1}g\sub{2}\,H \,=\, \,g\sub{1}\overline{g\sub{2}}\,H\,
$
for any $g\sub{1}\in G$.
Therefore  
$\overline{g\sub{1}\overline{g\sub{2}}} 
	\, = \,
\overline{g\sub{1}g\sub{2}}\,$ since 
$\overline{\tau}$ is well defined.
	\hfill$\Box$

\begin{definition}\label{def  frobenius} 
\emph{Given a  set of representatives  $T$ of $G/H$, the 
\emph{Frobenius representation} is the  representation}
\begin{equation}\label{equ  frobenius}
\gamma:
	\mar{G}{ \mathrm{S}\sub{\,T}}
\quad \mathrm{definied ~as}\quad
\every{g}{G}\every{t}{T}
	\big[~~
	{g\,} \longmapsto (~g^\gamma :  {t} \longmapsto\overline{g\,t}~)
	~~\big]\,.
\end{equation}
\end{definition}


To see that $\gamma$ is a representation, let $T$ be a
 set of representatives of
$G/H$ and  $g\sub{1}$ and $g\sub{2}$ be elements of $G$. Then
\vspace*{-5pt}
\[
(\,g\sub{1}g\sub{2}\,)^\gamma\,(t) 
	  \,=\, 
		 \overline{\,g\sub{1}g\sub{2}\,t} 
		\phantom{;}=\,
			\overline{ g\sub{1} \overline{g\sub{2}\,t}} 
	\, =\, \,\overline{\, g\sub{1} g\sub{2}^{\gamma}\,(t)  } 
	 \, =\, g\sub{1}^{\gamma} \circ g\sub{2}^{\gamma}\,(  \, t \,)\,,
\]
which follows from lemma \ref{lem  representative calculus}. Therefore
$
(\,g\sub{1}g\sub{2}\,)^\gamma\,=\, g\sub{1}^{\gamma}\circ g\sub{2}^{\gamma}\,.
$
Last, if $t$ is in $T$, then 
\,$\gamma: \marto{1}{1^\gamma}$  satisfies
$ 
1^\gamma \,(t) \, = \,\overline{1t} \, = \, \overline{t} \, = \, t
$
by definition, i.e.,  $1^\gamma\,=\, 1$ is in $\mathrm{S}\sub{\,T}$.
cf. \cite{Schupp}.

\section{The Diffraction and \emph{T}-Fibration Maps}

The Frobenius representation is fruitful in many branches of group theory, 
yet for us, it will give us the pivotal tools used to derive the 
mechanisms essential for our rewriting process.
To exploit it, we begin with the standard  representative system $T$
of $G/H$ for $G$. 
Now any group $G$ is bijective to $T\times H$:
\vspace{-5pt}
\[
G
	\,=\,
		\bigsqcup_{t\;\in\;T } tH
		\;\simeq_{\,_\Set}
		\bigsqcup_{t\;\in\;T }   \{\,t\,\}\times H
	\,=\,
		T\times H\,.
\vspace{-5pt}
\]
Since,  any $g\in G$, has a unique $t\in T$ so that 
$g\in tH$.  And, if  $g\,=\,t\,h$ and $\,t\,h =  t\,h'$ for 
$h,\,h'\in H$, then $h' = h$.
i.e., the map $\marto{g} {\ket{t,h} }$ is a well-defined 
bijection. For additional details see \cite{Rotman 1}.

\vspace{-5pt}

\begin{definition}\label{def  diffractor}
\emph{Given a set of representatives $T$ of $G/H$, the \emph{diffraction map} 
of $G$ by $T$ is the set-isomorphism}
\vspace{-5pt}
\begin{equation}\label{eq  nabla set iso}
\nabla
\,:\, \mar{G} {T\times H} 
\quad \mathrm{is ~defined~ by}\quad
	\nabla \,:\, \marto{g} {\ket{t,h} }
\end{equation}
\end{definition}

\vspace{-5pt}

If we think of the above
bijection metaphorically,  any element $g$ 
from its source $G$ can be witnessed with a ``prism"  $T$  as having a unique 
``spectrum"  $\ket{t,h}$ from its ``diffraction" $T\times H$.
Essentially, $\nabla$ diffracts the group $G$ 
into the set $T\times H$. Viewing $G$ internally, through the bijection, 
any element $g\in G$ 
 can also be uniquely decomposed  as $\,g\,=\,t\mspace{1mu}h$,
where $\ket{t,h}\in T\times H$.

\vspace{3pt}

The set of representatives $T$ under the product of $G$ 
is a sub-semigroup of $G$, and it is 
not a subgroup in general; this is why our goal is to ``diffract" $G$ as  
$T\times H$ such that it fosters  a multiplicative group structure on itself.
Therefore, to achieve this objective,  
consider the action $\hat{\varrho }\in \mor{\mathsf{Set}}{G\times G}{G}$
restricted to $T$,
induced by Cayley's representation $\varrho$ (cf. Definition 
 \ref{def  Cayley's representation}). Then,   rewrite the action  as 
\vspace{-10pt}
 \begin{alignat} {2}\label{eq  cayley mod frobinius}
\hat{\varrho}
	\;:\;
		\marto{\ket{g,t}}
			{
			} 
	&   
	g\mspace{1mu}t 
		\,=:\,
			\overline{g\mspace{1mu}t}\:
		 	\inv{(\,\overline{g\mspace{1mu}t}\,)}g\mspace{1mu}t
			& &
\end{alignat}
for any of $g\in G$ and $t\in T$. 
In the left-hand side  of  assignment (\ref{eq  cayley mod frobinius}) we have
rewritten the effect of Cayley's action over the pair 
$\ket{g,t}\in G\times T$ using the Frobenius representation.
Notice that $\overline{g \mspace{1mu}t} \in T$ by definition, and 
$\inv{(\,\overline{g\mspace{1mu}t}\,)}g\mspace{1mu}t \in H$ because 
for any  $t\in T$ and $g\in G$
\[
\overline{1}
 	\,=\,
		\overline{\,
			\inv{(\,\overline{g\,t}\,)\,}\,(\,\overline{g\,t}\,)
		}
	\,=\,
		\overline{\inv{(\,\overline{g\,t}\,)}\,g\,t\,}
\quad 
\Longrightarrow
\quad
	H
		\,=\,  
			\overline{1}\,H
		\,=\,
			\overline{\inv{(\,\overline{g\,t}\,)\,}\,g\,t\,}H
		\,=\,	
			\inv{(\,\overline{g\mspace{1mu}t}\,)}g\mspace{1mu}t\;H\,,
\vspace{-5pt}
\]
which follows from lemma \ref{lem  representative calculus}. 
Hence, assignment (\ref{eq cayley mod frobinius}) 
renders  the following tool:

\vspace{-3pt}

\begin{definition}
\emph{A \emph{$T$-fibration} of $G$, for 
a representatives $T$ of $G/H$, is the map}
\vspace{-5pt}
\[
\delta : \mar{   G\times T } {H}
\quad\textrm{defined  by}\quad
	\every{\ket{g,t}}{G\times T}
		\big[~
\,\delta(g,t)
	\, := \,
	\inv{(\,\overline{g\,t}\,)}g\,t
		~\big] \,.	
\] 
\end{definition}

A $T$-fibration $\delta$ is an abstract generalization of a fiber
bundle, at
its core. To see an analogous construction of $\delta$, consult 
\cite{Baumslag}. For us,  $\delta$ will allow 
a canonical description of the diffraction set-isomorphism 
$\nabla:\mar{G}{T\times H}$.

\vspace{-3pt}

\begin{Theorem}\label{the  nabla description}
If $T$ is a transversal of $G/H$ and $\delta:\mar{G\times T}{H}$ its
$T$-fibration, then the diffraction map $\nabla:\mar{G} {T\times H} $ 
uniquely decomposes into the pair 
\[
\nabla\,=\,\ket{\,\bar{\tau}\,,\,\delta\zero\,} 
\quad \text{such that}\quad 
	\every{g}{G}
		\big[~~
			\ket{
				\,\bar{\tau}\,,\,\delta\zero
			\,}\;( g )
				\,=\,
					\ket{\,
						\bar{g}\,,
							\inv{(\,\overline{g }\,)}g
					\,}
		~~\big] \,,
\]
where  
$\delta\zero$ is  $\, \delta(\aom,1):\mar{G}{H}$, and
$\bar{\tau}$ is the representative map  
\emph{(cf. definition \ref{def  representatives}).} 
\end{Theorem}


\noindent\emph{Proof}: Let $\delta : \mar{   G\times T } {H}$ be 
the $T$-fibration of $G/H$ and let $\nabla \,:\, \mar{G} {T\times H}$
be the diffraction map of $G$ by $T$. If $T$ is a transversal with $1\in T$, 
Cayley's action for any $g\in G$ on the identity   is just 
$\hat{\varrho}(g,1) = g$. Then
\[
\nabla(g)
	\,=\,
		\nabla \circ \hat{\varrho}\;(g,1)
	\,=\,
		\ket{\bar{g},\inv{(\,\overline{g }\,)}g}
	\,=\,
		\ket{\,\bar{\tau}(g) ,\delta (g,1)\,}
	\,=\,
		\ket{\bar{\tau},\delta\zero } \,(g)\,,
\]
where $\delta\zero$ is $\delta(\aom,1):\mar{G}{H}$. Now
given the projection maps 
$p\sub{H}:\mar{T\times H} {H}$ and 
$p\sub{T}:\mar{T\times H} {T}$, 
the pair
$\ket{\,\bar{\tau},\delta\zero}$ satisfies the commutative diagram
\vspace{-5pt}
\begin{equation*}\label{dia  factoring product}
\begin{aligned}[c]
\begin{tikzpicture}[description/.style={fill=white}]
\matrix (m) [matrix of math nodes, row sep=3.2em, column sep=2.5em, 
			text height=1.5ex, text depth=0.25ex] 
	{   & G  		  &         \\
	  T &  T\times H  & H  \,.      \\}; 
\path[->,font=\scriptsize]
	(m-1-2) edge node[above left] {$ \bar{\tau} $} (m-2-1)
	(m-1-2) edge node[above right ] {$ \delta\zero $} (m-2-3)
	 (m-2-2)edge node[above ]  {$ \phantom{asdf}p\sub{T} $} (m-2-1)
 	(m-2-2)edge node[above ]  {$ p\sub{H}\phantom{as} $} (m-2-3);
\path[->,dotted,font=\scriptsize]
 	(m-1-2)edge node[description]{$\ket{\,\bar{\tau},\delta\zero}$} (m-2-2);
\end{tikzpicture}
\end{aligned}
\end{equation*}
Therefore $\nabla$ is uniquely decomposable into the pair
$\,\ket{\,\bar{\tau},\delta\zero}$ by  uniqueness
of the universal property of direct products. 	
	\hfill$\Box$

\section{The Diffraction Representation}

In this penultimate section, we gather the tools we have built to describe
 a \emph{faithful} representation from $G$
into $\f{S}\sub{\,T\times H}$, in exchange for Cayley's representation.
Then, in the following last section, 
this representation will enable us to rewrite
the product of $G$ as a product in its 
canonical set $T\times H$ for a transversal $T$ of $G/H$.
We do this by first examining the general natural isomorphism
\vspace{-5pt}
\[
\mor{\cat{set}\,}{G\times T}{  H \,} 
	\, \cong \,
\mor{\cat{Grp}\,}{G}{ H^{\,T }\, }\,,
\vspace{-2pt}
\]
where $H^{\,T }:= \mor{\cat{set}\,}{T}{  H \,}$ is the group of maps from 
$G$ to $H^{\,T }$  under  functional point-wise multiplication, i.e.
\begin{equation}\label{eq  func prod grp}
	\every{f_1, \,{f_2}\,}{ H^{T }}
		\big[~f_1\cdot  {f_2\; }
	\,	: = \,
		\{~
			\ket{\,t\,,f_1(t)\, {f_2(t)}\;}
				~~|~ t\in T
		~\}
		~]\,.
\end{equation}
Notice the dual  
$\hat{\delta }\in \mor{\cat{set}\,}{G}{ H^{\,T }\, }$ of 
 a $T$-fibration $\delta:\mar{G\times T}{H}$, from the set of maps
 $\mor{\cat{set}\,}{G\times T}{  H \,} $, is described as 
\vspace{-2pt}
\begin{equation*}\label{def  dual t-fib}
\hat{\delta }
	:\mar{G} {\: H^{ \,T } }
\quad \text{such that}\quad  
\every{g}{G} 
\every{t}{T}
	~\Big[~~
\hat{\delta}:\marto{g}
	{ 
		\big(~
			\delta_g  
				: ~\marto{t}{ \delta(g,t) }
		~\big)
	}
	~~\Big]\,,
	\vspace{-5pt}
\end{equation*}
where $\delta_g$ is the map $\, \delta(g,\aom):\mar{T}{H}$ in $H^{\,T }$. 
Moreover,
the group structure of the image 
$
\hat{\delta}(G)
	\,:=\,
		\cset{\delta_g }{ g\in G },
$
of the dual   of  $\delta$ acts on the partition $T\times H$ of $G$. 
In general, any  
$f\in H^{\,T}$ can act as a permutation on the cartesian set
$T  \times H$, e.g.
\begin{equation}\label{eq  delta action}
\every{\ket{t,h} }{T\times H} 
	~\Big[~ 
		f\curvearrowright\ket{t,h}\,:=\, {\ket{\,t, \,f(t)\,h\,}}
	~\Big]\,.
\end{equation}

\begin{proposition}\label{lem  fibrated representation}
Given  the group  
$
\hat{\delta}(G),
$
from the image of the dual $\hat{\delta}$ of a $T$-fibration $\delta$, 
the assignment $\beta	:\mar{\hat{\delta}(G)}{\mathrm{S}\sub{\,T \times H}}$
defined as
\begin{equation}\label{eq  lem t-rep}
\every{\delta_g}{\,\hat{\delta}(G)} 
\every{\ket{t,h} }{\,T\times H} ~
\Big[~
\begin{tikzcd}[column sep=26]
		\delta_{g}   
			 \arrow[ |->,"\beta" ]{r}  & 	
		\Big(\,(\delta_{g})^\beta 
\end{tikzcd}
			{
				:
					\marto
						{\ket{t,h} }
							{
							\ket{t,\,\delta(g,t)\,h\,}
							}
				\,\Big)
			}
~~\Big]
\end{equation}
is an injective permutation representation of $\hat{\delta}(G)$ into 
$\f{S}\sub{\, T\times H }$.
\end{proposition}

\noindent \emph{Proof.} The assignment 
$
\beta:\mar{\hat{\delta}(G)}{\mathrm{S}\sub{\,T \times H}}
$
is a well-defined map by way of $\delta$. Now let
for any $\ket{t,h}\in T\times H$, by \emph{a priori} 
$
(\delta\sub{1})^\beta \ket{t,h} 
	\,=\,
		\ket{t,\delta\sub{1}(t)\, h}
	\,=\, 
		\ket{t,h}.
$
 Then 
$\beta$ is a homomorphism since 
\vspace{-3pt}
\[
\big(\delta_{g\sub{1}}\mspace{-6mu}\cdot \delta_{g\sub{2}}\big)^\beta \,
	\ket{t,h}
	\,=\,
	\ket{\,t,\delta(g\sub{1},t)\mspace{1mu}
	\delta(g\sub{2},t)\mspace{2mu}h }
	\,=\,
	\big(\delta_{g\sub{1} }\big)^\beta
		\ket{\,t, \delta(g\sub{2},t)\mspace{2mu}h\,}
	\,=\,
	\big(\delta_{g\sub{1}}\big)^\beta\!\circ\big(\delta_{g\sub{2}}\big)^\beta 
		\,\ket{t,h}		
\]
for any $g\sub{1},g\sub{2}\in G$.
To see that $\beta$ is injective, let 
${\delta_{g}}\in \kernel \beta$. Then
\begin{equation}\label{eq  lem t-rep ker}
(\delta_{g})^\beta \ket{t,h}
	\,=\,
		\ket{t, h}
	\,=\,
		\ket{t, \delta({g},t) h}\,.
\end{equation}
Equation
(\ref{eq  lem t-rep ker}) above is satisfiable if and only if 
$\delta_{g}(t) = 1$, 
i.e., $g = 1$. Therefore $\beta$ is an injective permutation representation 
of $\hat{\delta}(G)$ into  $\f{S}\sub{\, T\times H }$.
\hfill $\Box$

\vspace{10pt}

From proposition \ref{lem  fibrated representation} it follows that
 the group structure of $\hat{\delta}(G)$
 acts on the decomposition $T\times H$ of $G$. This result, along with the 
assistance of
the Frobenius representation, grants us our sought out extension 
of Cayley's representation 
to our permutation representation of $G$ into  $\f{S}\sub{\,T\times H}$.

\begin{definition} 
\label{def  diffracted representation}
\emph{The  \emph{diffraction representation} of $G$, 
by a set of representatives $T$ of $G/H$, is the assignment
$
\alpha :	\mar{G}{ \f{S}\sub{\,T\times H} }
$
defined by}
\vspace{-5pt}
\[
\every{g}{G} 
\!\every{\ket{t,h} }{T\times H} 
\Big[
\begin{tikzcd}[column sep=28]
		g   
			 \arrow[ |->,"\alpha" ]{r}  & 	
		\Big(\;
		(g^{\gamma}\times \id\sub{H})\,\circ \,(\delta_g)^\beta
\end{tikzcd}
			{
				:
					\marto
						{\ket{t,h} }
							{
							\ket{\,g^{\gamma}(t)\,,\,\delta(g,t)\,h\,}
							}
			\;\Big)
			}
~\Big]\,,
\]
\emph{where $\beta$ is the permutation representation described in 
proposition \ref{lem  fibrated representation}. } 
\end{definition}

\begin{lemma}\label{lem  T-homotopy calculus}
A $T$-fibration $\delta$, for a set of representatives $T$ of $G/H$,
satisfies 
\vspace{-4pt}
\[
\every{t}{T} 
\every{g\sub{1},g\sub{2}}{G} 
	~\Big[~~
		\,\delta( g\sub{1} g\sub{2}\,,t)
			\, = \, 
\delta\big(g\sub{1} ,\,g\sub{2}^{\,\gamma}\,(t)\big) \, \delta(g\sub{2},t)
	~~\Big]
\]
\end{lemma}

\noindent\emph{Proof}: If $\delta : \mar{   G\times T } {H}$ is 
the $T$-fibration associated to the representative $T$, then 
\vspace{-4pt}
\begin{equation*}
\delta\big(g\sub{1} ,\,g\sub{2}^{\,\gamma}\,(t)\big)
			\, \delta(g\sub{2},t)\, 
	\,=\,
		\inv{(\,\overline{g\sub{1}\,\overline{g\sub{2}\,t} }\,)}\,
		 \!g\sub{1}\,\overline{g\sub{2}\,t}\;
		 \,\inv{(\,\overline{g\sub{2}\,t}\,)}\, g\sub{2}\,t\,
		\, =\, 
		\inv{(\,\overline{g\sub{1} g\sub{2}\,t  }\,)}\,\,
		 g\sub{1}g\sub{2}\,t\,
	 \, =\, 
		\delta( g\sub{1}g\sub{2},t)\,
\end{equation*}
for any  $g\sub{1},g\sub{2}\in G$ and $t\in T$, since 
$\overline{g\mspace{1mu}t}\,: =\,g^\gamma(t)$.
\hfill$\Box$

\begin{Theorem}\label{the  fibrated rep}
The diffraction representation 
$\alpha :	\mar{G}{ \mathrm{S}\sub{ \,T\times H } } $
is an injective permutation representation  of $G$ into 
$\f{S}\sub{\, T\times H }$.
\end{Theorem}


\noindent\emph{Proof}: Let   $g\in G$ and  $\ket{t,h}\in T\times H$.
Then
$(g^{ \gamma }\times \idi{H} )\in \f{S}\sub{\,T\times H}$ since it acts as
\[
\begin{tikzcd}[column sep=48]
		\ket{t,h}  
			 \arrow[ |->,"g^{ \gamma}\times \idi{H} " ]{r}  & 	
				\ket{g^{ \gamma}(t), h} \,,
\end{tikzcd}
\]
which is permutation induced by the Frobenius representation $\gamma$. 
From Proposition \ref{lem  fibrated representation}, the map
$\delta_g \in \hat{\delta}(G)$ 
acts on the diffraction $T\times H$ as
$
(\delta_g)^\beta: \marto{\ket{t,h}}
			{\ket{\,t, \,\delta(t,g)\,h}}\,.
$
Thus 
$
(g^{\gamma}\times\idi{H})\,\circ \,(\delta_g)^\beta \in\f{S}\sub{T\times H }
$
is a permutation on the set $T\times H$ described as
\[
(g^{\gamma}\times \idi{H})\,\circ \,(\delta_g)^\beta :
					\marto
						{\ket{t,h} }
							{
							\ket{\,g^{\gamma}(t)\,,\,\delta(g,t)\,h\,}\,.
							}
\]
Now the action of the identity  $1\in G$ induced by $\alpha$  on 
$ \ket{h,t} $  is given as
\vspace*{-5pt}
\begin{alignat} {2}
\alpha(1)\,\ket{t,h} \, := \, 
  (1^{\gamma}\times \idi{H}) \circ (\delta\sub{1})^\beta \;\ket{t,h} \,=\,
		\ket{\,1^\gamma(t)\,,\,\delta(1,t)\,h\,}
	 \, =\, 
		\ket{\,\overline{t}\,,\,\inv{(\,\overline{ t}\,)} \,t\,h\,}
	&  \, =\, 
		\ket{ 
			t ,h
		 }
			& & 
			\notag
\end{alignat}
since $\bar{t} = t$. Therefore $\alpha(1) \, = \, 1$ in
$\f{S}\sub{T\times H}$.
Moreover, given $g\sub{1},g\sub{2}\in G$,  the action of $g\sub{2}$ 
followed by the action of $g\sub{1}$ on $ \ket{t,h} $ by way of $\alpha$ is
\vspace{-3pt}   
\begin{alignat} {2}
\alpha(g\sub{1})\circ \alpha(g\sub{2}) \;\ket{t,h}
&\,=\,
\big(\,
	 (g\sub{1}^{\gamma}\times \idi{H})\circ (\delta_{g\sub{1}})^\beta
\,\big)
	\circ 
		\big(\,
			(g\sub{2}^{\gamma}\times \idi{H})\circ(\delta_{g\sub{2}})^\beta
		\,\big)\;
			\ket{t,h}
			\notag\\ 
	&\,=\,\big(
			(g\sub{1}^{\gamma}\times \idi{H})\,\circ(\delta_{g\sub{1}})^\beta
	 \,\big)\; 
		\ket{\,g\sub{1}^{\,\gamma}\,(t)\,,\,\delta(g\sub{2},t)\,h\,}
		\notag\\
	& \, =\, 
		\ket{\,
			g\sub{1}^{\,\gamma}\circ g\sub{2}^{\,\gamma}\,(t)
			\,,\,\delta( g\sub{1} \,,\,g\sub{2}^{\,\gamma}\,(t)\,)
			\, \delta(g\sub{2},t)\,h
		\,}\,.
			& &
			\notag
\end{alignat}
But  
$
g\sub{1}^{\,\gamma}\circ g\sub{2}^{\,\gamma} \;(t)
	\,=\,
		(\, g\sub{1}\,g\sub{2}\,)^{\,\gamma}\,(t)
$\,
and
$
\delta\big(g\sub{1} \,,\,g\sub{2}^{\,\gamma}\,(t)\big)
			\, \delta(g\sub{2},t)\, 
	 \, =\, 
		\delta( g\sub{1} \,g\sub{2}\,,t)\,
$
by the Frobenius representation $\gamma$ and  lemma
\ref{lem  T-homotopy calculus}, respectively. Therefore,
putting these together, 
we get 
\vspace{-3pt}
\begin{alignat} {2}
\alpha(g\sub{1})\circ \alpha(g\sub{2}) \;\ket{t,h}
	& \, =\, 
		\ket{\,( g\sub{1}g\sub{2})^{\gamma}\,(t)
			\,,\,
			\delta(\,g\sub{1} g\sub{2},t)\;h
		\,}
	 \, =\, 
		\big(\,( g\sub{1}g\sub{2})^{\gamma}\times \idi{H}\big)
		\circ 
		(\delta_{\,g\sub{1}  g\sub{2}})^\beta
		  \; \ket{t,h}
			& &
			\notag\\
	& \, =\, 
		\alpha(\,g\sub{1} g\sub{2}\,)\; \ket{t,h}\,,
			& &
			\notag
\end{alignat}
i.e.,
$\alpha(\,g\sub{1}g\sub{2}\,)
	\,=\,
		\alpha (g\sub{1}) \circ \alpha(g\sub{2})\,.
$
Hence $\alpha$ is a representation of 
${G}$ into $ \f{S}\sub{\,T\times H}$.
Last, $\alpha$ is injective if 
$
\ket{\,g^{\gamma}(t)\,,\,\delta_g(t)\,h\,}
	\,=\,
		\ket{t, h}\,,
$
i.e., 
$\overline{g\mspace{1mu}t}=t$ and $\delta_{g}(t) = 1$
if and only if $g = 1$. 
Therefore, $\alpha$ is   injective as well.
	\hfill $\Box$

\vspace{5pt}


\section{Rewriting in a Diffracted Group}
\label{sec  dif gp}

Given a transversal set $T$ of $G/H$ for a group $G$,
in this last section, we describe the group structure for the set $T\times H$ 
that emerges from the bijection of the diffraction  map
$\nabla:\mar{G}{T\times H}$ and its diffraction representation 
$\alpha :	\mar{G}{ \mathrm{S}\sub{ \,T\times H } } $.  Moreover,
 we will also  show
that $\nabla$ can be lifted to a universal group-isomorphism
$
G\, \simeq\sub{\cat{Grp}} T\times H\,.
$ 
We begin by applying a group-action isomorphism on sets, cf. \cite{Mac Lane}.

\begin{definition}\label{def  action cat}
\emph{Given a group $G$, the category 
$G\textrm{-\textsf{Set}}$ of (left) \emph{group actions of $G$ on 
$\mathsf{Set}$} consists of 
the following objects and morphisms:}  

\vspace*{5pt}

\noindent $\bullet~$ 
\emph{Objects in  $G\textrm{-\textsf{Set}}$ are pairs
$\ket{\hat{\varrho},X}$, where $X\in \textsf{Set}$ and 
$
\hat{\varrho} \in \mor{\mathsf{Set}} {G\times X}{X}
$
satisfy}
\begin{description}
\item[\cat{i\,.\phantom{i}}] 
$
\every{x}{X} 
	~\big[~~
		\hat{\varrho} (\mspace{2mu}1,x\mspace{1mu}) = x
	~~\big]
$

\item[\cat{ii\,.}]
$
\every{g\sub{1},\,g\sub{2}}{G} \every{x}{X} 
	~\big[~~ 
		\hat{\varrho} (\,g\sub{1}\mspace{1mu}g\sub{2},x\,) 
		\,=\, 
		\hat{\varrho} 
		\big(\,g\sub{1},\mspace{1mu}
		\hat{\varrho} (\mspace{1mu}g\sub{2},x)\,\big) 
	~~\big]
$

\end{description}
\emph{We refer to  objects $\ket{\hat{\varrho} ,X}$ from 
$G\textrm{-\textsf{Set}}$ as (left) \emph{actions} of $G$ on sets.}

\vspace*{5pt}

\noindent $\bullet~$
\emph{Morphisms in $G\textrm{-\textsf{Set}}$  are arrows between any 
two actions
$
\nabla : \mar{ \ket{\hat{\varrho} ,X}}  { \ket{\hat{\alpha} ,Y} }
$
such that $\nabla \in \mor{\mathsf{Set} }{X}{Y}$ and 
the following diagram commutes:}

\[
\begin{tikzpicture}[description/.style={fill=white}]
\matrix (m) [matrix of math nodes, row sep=1.8em, column sep=3em, 
			text height=1.5ex, text depth=0.25ex] 
	{ G \times X & X  			\\
	  G \times Y & Y		\\ }; 
\path[->,font=\scriptsize]
	(m-1-1) edge node[left] {$ \id \times \nabla $} (m-2-1)
	(m-1-1) edge node[above] {$ \hat{\varrho}  $} (m-1-2)
 	(m-1-2)edge node[right] {$ \nabla $} (m-2-2)
    (m-2-1)edge node[above] {$ \hat{\alpha} $} (m-2-2);
\end{tikzpicture}
\] 
\vspace{-20pt}

\noindent\emph{A morphism $\nabla$ in
$G$-$\mathsf{Set}$ is said to be a \emph{$G$-$\mathsf{Set}$ isomorphism} 
if $\nabla$ is also a bijective map.} 
\end{definition}

\begin{Theorem}\label{the G-set isomorphism}
The diffraction map $\nabla:\mar{G}{T\times H}$ 
is a $G$-$\Set$ isomorphism between the actions induced by Cayley's 
representation 
$
\varrho  \,: \mar{G }{S\sub{G} }
$
and the diffracted representation 
$
\alpha \,:  \mar{ G }{S\sub{T\times H} }\,,
$
for a transversal $T$ of $G/H$.
\end{Theorem}

\noindent \emph{Proof}: Let 
$
\hat{\varrho}  \,: \mar{G\times G }{ G }
$
and
$
\hat{\alpha} \,: \mar{ G\times (T\times H)}{T\times H}
$
be the actions induced by Cayley's representation $\varrho$ and the
diffraction representation $\alpha$ from Theorem \ref{the  fibrated rep},
respectively. 
Then our bijection  $\nabla : \mar{G} {T\times H}$ satisfies the commutative
diagrams
\begin{equation*}\label{dia  G set isomorphism} 
\begin{aligned}[c]
\begin{tikzpicture}[description/.style={fill=white}]
\matrix (m) [matrix of math nodes, row sep=1.8em, column sep=2.5em, 
			text height=1.5ex, text depth=0.25ex] 
	{ G\times G  			  & G				\\
	  G\times\,(T\!\times\!H) & (T\!\times\!H)	\\ 
	}; 
\path[->,font=\scriptsize]
	(m-1-1) edge node[left] {$ \idi{G} \times \nabla~ $} (m-2-1)
	(m-1-1) edge node[above] {$ \hat{\varrho} $} (m-1-2)
 	(m-1-2)edge node[right] {$ \nabla $} (m-2-2)
    (m-2-1)edge node[auto] {$ \hat{\alpha} $} (m-2-2);
\end{tikzpicture}
\end{aligned}
\quad \text{since}\quad
\begin{aligned}[c]
\begin{tikzpicture}[description/.style={fill=white}]
\matrix (m) [matrix of math nodes, row sep=1.8em, column sep=2em, 
			text height=1.5ex, text depth=0.25ex] 
	{ \ket{g,\,k\,}  			  & 
			\overline{gt}\;
			\inv{(\,\overline{g\,t}\,)}\;g\,t\,h
													\\
	  \ket{ g\,,\ket{t,h}\,}  & 
	  		\ket{
			\overline{gt}\,,\,
			\inv{ (\,\overline{g\,t}\,)} g\,t\,h
			}										\\ 
	}; 
\path[|->,font=\scriptsize]
	(m-1-1) edge node[left] {$ \idi{G} \times \nabla~ $} (m-2-1)
	(m-1-1) edge node[above] {$ \hat{\varrho} $} (m-1-2)
 	(m-1-2)edge node[right] {$ \nabla $} (m-2-2)
    (m-2-1)edge node[auto] {$ \hat{\alpha} $} (m-2-2);
\end{tikzpicture}
\end{aligned}
\end{equation*}
for any $g,k\in G$ and  $k=t\,h$, where $\ket{t,h}\in T\times H$. i.e., 
$\nabla$ is a $G$-$\Set$ isomorphism.
\hfill$\Box$

\vspace{10pt}

Theorem \ref{the G-set isomorphism} not only stipulates that  
$G$ and $T\times H$ are $G$-$\Set$ isomorphic, but that the
diffracted representation $\alpha$ expresses  Cayley's representation 
when $G$  decomposes as $\nabla(G)\,=T\times H$ as a set. 
Fortunately, this relationship is more than happenstance; the diffraction  
$T\times H$ of $G$ borrows the group structure of $G$ through  
the diffraction representation $\alpha$ since  Cayley's representation
$\varrho$ exposes the product structure 
of  $G$.  As a consequence, following the commutativity of the
diagram in (\ref{dia  Cayley Product}) from the introduction,
 the product ``$\bullet$" in $G$ can be reconstructed  as 
\begin{equation}\label{dia  difracted Product}
\begin{aligned}[c]
\begin{tikzpicture}[description/.style={fill=white}]
\matrix (m) [matrix of math nodes, row sep=1.8em, column sep=2.5em, 
			text height=1.5ex, text depth=0.25ex] 
	{ G \times G & G  \\
	G \times(T\times H) & T\times H\\ }; 
\path[->,font=\scriptsize]
	(m-1-1) edge node[left]  {$  \id \times \nabla $} (m-2-1)
	(m-1-1) edge node[above] {$ `` \bullet " $} 	(m-1-2)
    (m-2-1)edge node[auto]   {$ \hat{\alpha} $} (m-2-2);
    \path[<-,font=\scriptsize]
    (m-1-2)edge node[right]  {$ \nabla\,\inverse  $} (m-2-2);
\end{tikzpicture}
\end{aligned} 
\end{equation}

\begin{Theorem}\label{the  diffracted group}
If $\delta$ is the $T$-fibration associated to  Frobenius representation
$\gamma$ for a transversal $T$ of $G/H$, then 
$T \times H$ becomes a group under the product defined by
\vspace{-5pt}
\begin{equation}\label{eq  diffration group operation}
\every{ \ket{t\sub{1},h\sub{1} },\,\ket{t\sub{2},h\sub{2}} } {T\times H} 
	~\Big[~~
\ket{t\sub{1},h\sub{1} }\by \ket{t\sub{2},h\sub{2} }	
	\: :=\:
		\ket{\,
			t\sub{1}^{\gamma}\!\circ\! \,h\sub{1}^{\gamma} \,(t\sub{2})\,,
			\,\delta(t\sub{1}h\sub{1},t\sub{2})\,h\sub{2}
		\,} 
		~~\Big]
\end{equation} 
\end{Theorem}

\noindent\emph{Proof}. Given the diffraction map 
$\nabla : \mar{G} {T\times H}$
of $G$ with respect to a transversal $T$ 
and the diffraction representation 
$
\alpha \,:  \mar{ G }{T\times H}\,
$ 
of $G$ by $T$, it follows from diagram (\ref{dia  difracted Product})
by mean of theorem \ref{the G-set isomorphism}
that the following diagram commutes:

\vspace{-5pt}
\begin{equation*}\label{dia  proof difracted Product}
\begin{aligned}[c]
\begin{tikzpicture}[description/.style={fill=white}]
\matrix (m) [matrix of math nodes, row sep=1.6em, column sep=.7em, 
			text height=1.5ex, text depth=0.25ex] 
	{   
	G\times G  & 			   & & &           & 	G   \\
			   & (T\times H)^2 & & & T\times H &			 \\
	G\times(T\times H) 
			   & 			   & & &           & T\times H\\ }; 
\path[->,font=\scriptsize]
	(m-1-1) edge node[right]  {$\id\times\nabla $} (m-3-1)
	(m-1-6) edge node[above left]  {$\nabla $} (m-2-5)
	(m-1-6) edge node[above left]  {$\nabla $} (m-2-5)
	(m-2-2) edge node[above]   {$``\cdot" $} (m-2-5)
	(m-2-2.160) edge node[above right]
		{$\nabla\inverse \!\times\nabla\inverse$} (m-1-1)
	(m-1-1) edge node[above]{$  \mspace{38mu}``\bullet" $} 	(m-1-6)
	(m-3-6)edge node[left]  {$\nabla\inverse $} (m-1-6) 
    (m-3-1)edge node[below right]   {$ \hat{\alpha} $} (m-3-6);
\path[<-,densely dashed,font=\scriptsize]
	(m-3-1.15) edge node[ right]
	{$ \mspace{10mu}(\nabla\inverse\; )\times \id$} (m-2-2.200)
    (m-2-5)edge node[ left]  {$  \id \mspace{10mu}$} (m-3-6);
\end{tikzpicture}
\end{aligned} 
\end{equation*} 
\vspace{-7pt}

\noindent Moreover, letting  $g\sub{1} = t\sub{1}h\sub{1}$ and $g\sub{2} = t\sub{2}h\sub{2}$ in $G$ where $\ket{t\sub{1},h\sub{1}}, \ket{t\sub{2},h\sub{2}}\in T\times H$, then chasing the  central arrow $``\cdot"$ 
through $\hat{\alpha}$ describes the product in $T\times H$  from $G$ 
as
\[
\begin{tikzcd}[column sep=28]
		\big\langle\, \ket{t\sub{1},h\sub{1}},\ket{t\sub{2},h\sub{2}} 
		\,\big\rangle
		\arrow[ |->,"  ``\,\cdot\mathrm{\,"}  "  ]{r}  & 	
		\ket{\,
			t\sub{1}^{\gamma}\!\circ\! \,h\sub{1}^{\gamma} \,(t\sub{2})\,,
			\,\delta(t\sub{1}h\sub{1},t\sub{2})\,h\sub{2}
		\,} \,.
\end{tikzcd}
\]
Therefore the set $T\times H$ becomes   
a group under formulation (\ref{eq  diffration group operation})
since  $G$ is a group.
\hfill $\Box$

\begin{definition} \label{def  external diffracted group}
\emph{Let $ T$ be a transversal of $G/H$ for  $G$. 
The  \emph{diffracted group of $G$}  by $T$, denoted by
$
	 {T \,\triangledown H}\,,
$
is the group structure on $T\times H$ described by the product in 
formulation (\ref{eq  diffration group operation}). And, the 
product in
$
	 {T \,\triangledown H}
$
is termed the (\emph{external\,$)$ bequeath product}.  
}
\end{definition}

Theorem \ref{the  diffracted group} also allows us to use the Frobenius
 representation $\gamma$
and its $T$-fibration $\delta$ to rewrite the product of any two
elements $g\sub{1},g\sub{2} \in G$ in 
terms of a transversal $T$ and its associated subgroup $H$.
More precisely, given any two elements
 $g\sub{1} = t\sub{1}h\sub{1}$ and $g\sub{2} = t\sub{2}h\sub{2}$ in $G$
 of the form $\ket{t\sub{1},h\sub{1}}, \ket{t\sub{2},h\sub{2}}\in T\times H$, 
we can use the bequeathed rewriting process to express 
their internal product in terms of
elements in  $T$ and $H$ as
 \vspace{-3pt}
\[
g\sub{1}\mspace{1mu} g\sub{2} 
	=
			\,t\sub{1}^{\gamma}\!\circ\! \,h\sub{1}^{\gamma} \,(t\sub{2})
			\,\delta(t\sub{1}h\sub{1},t\sub{2})\,h\sub{2}
	=
			\,\overline{t\sub{1} h\sub{1} t\sub{2}}
			\,\delta(t\sub{1}h\sub{1},t\sub{2})\,h\sub{2}\,,
\]
which is  
$g\sub{1}\mspace{1mu} g\sub{2} = t\sub{1}h\sub{1}t\sub{2}h\sub{2}$;
nevertheless, elucidating
 $\overline{t_1 h_1 t_2}\in T$ and $\delta(t_1h_1,t_2)\,h_2 \in H$.

\begin{definition} \label{def  bequeathed product}
\emph{Let $ T$ be a transversal of $G/H$ for  $G$  and 
$g\sub{1},g\sub{2} \in G$ such that $g\sub{1} = t\sub{1}h\sub{1}$ and 
$g\sub{2} = t\sub{2}h\sub{2}$,
where $t\sub{1},t\sub{2}\in T$ and $h\sub{1},h\sub{2}\in H$. 
The (\emph{internal\,$)$ bequeathed product} of $g\sub{1}$ and $g\sub{2}$ is 
the rewriting of $g\sub{1}g\sub{2}$ 
in terms of elements from  $T$ 
and $H$ defined by}
\[
g\sub{1}\mspace{1mu} g\sub{2} 
	\,:=\,
\,\overline{t\sub{1} h\sub{1} t\sub{2}}
\,\delta(t\sub{1}h\sub{1},t\sub{2})\,h\sub{2}\,.
\]
\end{definition}

Furthermore, Theorem \ref{the  diffracted group} not only dispenses the 
essential framework to convert the set $T\times H$ into a group,
but it also assists in the  lifting of the  diffraction bijection 
$\nabla \,:\, \mar{G} {T\times H}$ to a group-isomorphism 
$\nabla:\mar{G}{{T\,\triangledown H}}$.

\begin{Theorem}\label{pro  diffractor isomorphism}
If $G$ is a group with a subgroup $H$ and an embedding 
$\iota\sub{1}:\mar{H}{G}$,  
then for any  transversal\, $T$ of\, $G/H$ with an embedding 
$\iota\sub{2}:\mar{H}{T\,\triangledown H}$,  
there is a unique isomorphism 
$\nabla:\mar{G}{{T\,\triangledown H}}$ such that 
$\nabla \circ\iota\sub{1} = \iota\sub{2}$, 
i.e., the  diagram commutes:
\vspace{-5pt}
\begin{equation*}\label{dia  proof univ diffraction}
\begin{aligned}[b]
\begin{tikzpicture}[description/.style={fill=white}]
\matrix (m) [matrix of math nodes, row sep=.3em, column sep=2em, 
			text height=1.5ex, text depth=0.25ex] 
	{    &  G  			  	           &  & G					   \\
	  H  & \phantom{a}\mspace{45mu}          &  & \mspace{30mu}\phantom{b}\\ 
	     & T\,\triangledown H      &  & T\times H                    
    \\}; 
\path[->,font=\scriptsize]
	(m-2-1) edge node[above] {$ \iota\sub{1} $} (m-1-2.190)
	(m-2-1) edge node[below] {$  \iota\sub{2} $} (m-3-2.170)
 	(m-1-2)edge node[right] {$\nabla $} (m-3-2);
 \path[<-|,dotted,font=\scriptsize]
    (m-1-2)edge (m-1-4)
    (m-3-2)edge  (m-3-4)
	(m-2-2.east) edge node[description] {$ \cat{Grp} $} (m-2-4.west);
\path[->,font=\scriptsize]
 	(m-1-4)edge node[left] {$\nabla$} (m-3-4);
\end{tikzpicture}
\end{aligned}
\end{equation*}
\end{Theorem}

\noindent\emph{Proof}: 
Let $T$  be a transversal of $G/H$, and $\nabla:\mar{G}{T\times H}$ be its
diffraction bijection.
Given the identity $1\in G$, then $\nabla(1) = \ket{1,1}$. Now 
let $g\sub{1},g\sub{2}$ be elements in $G$ such that each has a unique 
decomposition 
$
g\sub{1} \,=\, t\sub{1}\,h\sub{1}
$
and
$
g\sub{2} \,=\, t\sub{2}\,h\sub{2}\,,
$
where $t\sub{1},t\sub{2}\in T$ and $h\sub{1},h\sub{2}\in H$. Then 
 $\nabla$ extends to a homomorphism from $G$ into $T \,\triangledown H$ 
since
\begin{alignat}{2}
\nabla(g\sub{1}g\sub{2})
&\,=\;
	 \ket{\,\overline{t\sub{1}h\sub{1}\, t\sub{2}h\sub{2}} \,, 
	 \delta(t\sub{1}h\sub{1}\, t\sub{2}h\sub{2},1)\,}
&&\,=\;
	 \ket{\, 
	 \overline{t\sub{1}h\sub{1}\, t\sub{2} \overline{h\sub{2} }} \,,\, 
	 \delta(t\sub{1}h\sub{1} ,\, \overline{ t\sub{2}\overline{h\sub{2} } } )
	 				\,\delta( t\sub{2}h\sub{2}, 1)
	 \,}
			\notag \\
&\,=\,
	\ket{ \,
	 	t\sub{1}^{\;\gamma}\circ h\sub{1}^{\;\gamma} ( t\sub{2}) \,, 
	 		\delta(t\sub{1}h\sub{1}, \, t\sub{2} )\,h\sub{2}
	 \,}	
&&\,=\,
			\ket{t\sub{1},h\sub{1}}\by\ket{t\sub{2},h\sub{2} }		
			\notag	\\
&\,=\,
			\nabla(g\sub{1} )\by \nabla( g\sub{2})\,
			\notag	
\end{alignat}
by Lemma \ref{lem  representative calculus}, 
Lemma \ref{lem  T-homotopy calculus} and Theorem \ref{the  diffracted group}.
Moreover $\nabla$ is an isomorphism 
since any $g\in \kernel \nabla$ must satisfy $g=1$. 
The universal mapping property
$\nabla \circ\iota\sub{1} = \iota\sub{2}$ follows by composition,
where the existence and
uniqueness of   $\nabla$
  is provided by  Theorem 
\ref{the  nabla description}.
\hfill$\Box$


\medskip

\noindent
Mathematics Department, 
New York City College of Technology, Brooklyn
 NY 11201

\medskip

\noindent {\it e-mail \/}:
  nyzapata@gmail.com
\medskip


\end{document}